\documentclass[12pt]{amsart}
\usepackage{amssymb}
\usepackage{amscd}

\newtheorem{Theorem}{Theorem}[section]
\newtheorem{Lemma}[Theorem]{Lemma}
\newtheorem{Corollary}[Theorem]{Corollary}
\newtheorem{Proposition}[Theorem]{Proposition}
\theoremstyle{definition}

\newtheorem{Example}[Theorem]{Example}
\theoremstyle{remark}
\newtheorem{Remark}{Remark}

\font\sy=cmsy10

\font\ym=msbm10  
\newcommand{\cD}{{\hbox{\sy D}}}
\newcommand{\cE}{{\hbox{\sy E}}}

\newcommand{\cH}{{\hbox{\sy H}}}

\newcommand{\cK}{{\hbox{\sy K}}}
\newcommand{\cL}{{\hbox{\sy L}}}

\newcommand{\cO}{{\hbox{\sy O}}}

\newcommand{\cT}{{\hbox{\sy T}}}

\newcommand{\cV}{{\hbox{\sy V}}}

\newcommand{\C}{{\text{\ym C}}}

\newcommand{\N}{\text{\ym N}}

\newcommand{\R}{\text{\ym R}}

\newcommand{\End}{\hbox{\rm End}}
\newcommand{\Hom}{\hbox{\rm Hom}}

\usepackage[dvips]{graphicx}

\title[]
{Oriented Kauffman diagrams and universal quantum groups}
\author[Yamagami Shigeru]{Yamagami Shigeru}

\begin{document}
\maketitle   
\begin{center}
Department of Mathematics and Informatics 
\end{center}
\begin{center}
Ibaraki University 
\end{center}
\begin{center} 
Mito, 310-8512, JAPAN 
\end{center}    

\begin{abstract}
We study 
the tensor category of oriented Kauffman diagrams 
and determine fiber functors on them as well as 
the associated Hopf algebras. 
\end{abstract}
\bigskip
                           

\section{Introduction}
In our previous paper \cite{CDA, FTL}, 
we determined fiber functors 
on Temperley-Lieb categories with the help of universality 
property on rigidity and classification results 
on bilinear forms. The associated 
Hopf algebra as a consequence of Tannaka-Krein duality 
is identified with 
the algebraic quantum group of Dubois-Violette and Launer. 

When C*-structure is entailed, 
unitary fiber functors are also classified 
with the associated compact quantum groups 
isomorphic to the universal quantum groups of orthogonal type 
due to Wang \cite{D-W} and Banica \cite{Ba}. 
There is another related class 
of compact quantum groups, called 
the universal quantum group of unitary type, 
investigated by the same authors. 

Here we shall study an oriented version of Temperley-Lieb 
categories and obtain an analogous classification 
of fiber functors on them. 
Though we have failed in identifying the associated 
algebraic quantum groups, when restricted to the unitary case, 
they turn out to be the universal quantum groups of 
unitary type, thus revealing 
geometric structures behind them. 

\section{Oriented Kauffman Diagrams}
Let $D$ be a Kauffman diagram of type $(m,n)$ 
(\cite{Kau, Kau2}), i.e., 
$D$ is the isotopy class of planar strings in a rectangle 
with the strings having $m$ terminal points on 
the upper bounding line 
and $n$ terminal points on the lower bounding line.
Thus $D$ contains $(m+n)/2$ strings in total and 
there are $2^{(m+n)/2}$ possibilities in the choice of orientations 
of strings in $D$. 
A Kauffman diagram with a specific orientation is called 
an \textbf{oriented Kauffman diagram}. 
For an explicit description of orientation, 
we assign one of symbols $X$ and $X^*$ to each point 
so that $X$ (resp.~$X^*$) indicates starting 
(resp.~ending) for upper bounding points and 
ending (resp.~starting) for lower bounding points.
By this coding, an orientation for $D$ produces 
two words $(X_1,X_2,\dots, X_m)$ and  
$(X'_1,X'_2,\dots, X'_n)$ of two letters $\{ X, X^* \}$
corrsponding 
to upper and lower sequences of vertices. 
Clearly the set of such words is  
identified with the free product monoid $\N * \N$ 
($\N = \{ 1, 2, \dots\}$ 
being the additive monoid of natural numbers). 
For an element $w \in \N * \N$, we use the notation $X^w$ 
to stand for the corresponding word. 
Let $w, w' \in \N *\N$ be defined by 
\[
X^w = (X_1, \dots, X_m), 
\qquad 
X^{w'} = (X'_1, \dots, X'_n).
\]
The pair $(w,w')$ is called the type of 
an oriented Kauffman diagram. 
Given words $w, w' \in \N * \N$, let $K_{w,w'}$ be the set of 
oriented Kauffman diagrams of type $(w,w')$ and $\C[K_{w,w'}]$ be 
the free complex vector space generated by the set $K_{w,w'}$ 
($\C[\emptyset] = \{ 0\}$ by definition). 
Let $|w|_+$ (resp.~$|w|_-$)  
be the number of $X$'s (resp.~$X^*$'s) in $X^w$. 
Then the set $K_{w,w'}$ is non-empty if and only if both of 
$|w|_+ + |w'|_-$ and $|w|_- + |w'|_+$ are even numbers,    
which can be seen by an easy induction argument.

Just as in the unoriented case, 
we introduce a (strict) tensor category 
$\cO_{d_L,d_R}$ parametrized by $d_L, d_R \in \C^\times$: 
\begin{enumerate}
\item
Objects are exactly the symobols $\{ X^w\}$ with $w \in \N * \N$. 
\item
Hom-sets are set to be $\Hom(X^w,X^{w'}) = \C[K_{w,w'}]$ 
and the operation of composition 
$\Hom(X^{w'},X^{w''})\times \Hom(X^w,X^{w'}) 
\to \Hom(X^w,X^{w''})$ is defined by 
the concatenation of planar strings with clockwise 
(resp.~anticlockwise) loops replaced by $d_R$ (resp.~$d_L$).  
\item
The tensor product for morphisms is the linear extension of 
the horizontal juxtaposition of diagrams in $K_{w,w'}$. 
\end{enumerate}

Notice that the unit object $I$ of $\cO_{d_L,d_R}$ is associated 
to the empty word, i.e., $I = X^\emptyset$.
 
The tensor category obtained in this way is rigid as in the case 
of Temperley-Lieb categories and bears a kind of universality 
on rigidity.

Let $\epsilon_X: X\otimes X^* \to I$ and 
$\epsilon_{X^*}:X^*\otimes X \to I$ be basic arcs in 
$\cO_{d_L,d_R}$ with 
$\delta_X: I \to X^*\otimes X$ and $\delta_{X^*}:I \to X\otimes X^*$ 
the associated copairings. 

\begin{Lemma}
Any (strictly) monoidal functor $F$ of 
$\cO_{d_L,d_R}$ into a tensor category
$\cT$ is uniquely determined by the choice 
\[
F(\epsilon_X): F(X)\otimes F(X^*) \to I, 
\quad 
F(\epsilon_{X^*}): F(X^*)\otimes F(X) \to I
\]
of rigidity pairings which satisfy the identities
\[
F(\epsilon_X) F(\delta_{X^*}) = d_L 1_I, 
\quad
F(\epsilon_{X^*}) F(\delta_X) = d_R 1_I.
\]
Note here that $F(\delta_X)$ and $F(\delta_{X^*})$ are 
characterized as the rigidity coparings associated with 
$F(\epsilon_X)$ and $F(\epsilon_{X^*})$ respectively. 

Conversely, given rigidity pairings 
$\epsilon_Y: Y\otimes Y^* \to I$ and $\epsilon_{Y^*}: 
Y^*\otimes Y \to I$ in a tensor category $\cT$, satisfying 
the relations 
\[
\epsilon_Y\delta_{Y^*} = d_L 1_I, 
\quad 
\epsilon_{Y^*}\delta_Y = d_R 1_I,
\]
there exists a monoidal functor $F: \cO_{d_L,d_R} \to \cT$ 
such that $F(X) = Y$, $F(X^*) = Y^*$, 
$F(\epsilon_X) = \epsilon_Y$ and 
$F(\epsilon_{X^*}) = \epsilon_{Y^*}$.
\end{Lemma}

Although the category $\cO_{d_L,d_R}$ has apparently 
two parameters, 
one freedom of them is superficial as we shall see below.

\begin{Lemma}
Let $X^w$ ($w \in \N * \N$) be an object of $\cO_{d_L,d_R}$. 
Then 
$\dim \End(X^w) = 1$ if and only if $X^w = X^n$ or 
$X^w = (X^*)^n$ for some $n \in \N$. 
\end{Lemma}

\begin{proof}
In fact, $\End(X\otimes X^*)$ and $\End(X^*\otimes X)$ are 
two-dimensional, which are included in 
$\End(X^w)$ unless $X^w = X^n$ or $X^w = (X^*)^n$. 
\end{proof}

Let $F: \cO_{d_L,d_R} \to \cO_{d'_L,d'_R}$ be an equivalence 
of tensor categories, i.e., $F$ is an essentially surjective 
and fully faithful tensor functor. 
By the above lemma, we then have 
$F(X) = X^n$ or $F(X) = (X^*)^n$ for some integer $n \geq 1$. 
The case $n \geq 2$, however, contradicts with the essential 
surjectivity of $F$ because it implies $X \not\cong F(X^w)$ 
for any $w \in \N *\N$. 
Thus we have alternatives $F(X) = X$ or $F(X) = X^*$. 

The operation of $F$ on morphisms is then determined by 
the effect on basic arcs $\epsilon_X$ and $\epsilon_{X^*}$: 
\[
F(\epsilon_X) = \lambda \epsilon_X, \quad 
F(\epsilon_{X^*}) = \mu \epsilon_{X^*}
\quad 
(F(\delta_X) = \lambda^{-1}\delta_X, \quad 
F(\delta_{X^*}) = \mu^{-1} \delta_{X^*})
\]
or 
\[
F(\epsilon_X) = \lambda \epsilon_{X^*}, \quad 
F(\epsilon_{X^*}) = \mu \epsilon_X
\quad
(F(\delta_X) = \lambda^{-1}\delta_{X^*}, \quad 
F(\delta_{X^*}) = \mu^{-1} \delta_X). 
\]
From the obvious equality 
\[
F(\epsilon_X\delta_{X^*}) = \epsilon_X\delta_{X^*}, 
\quad 
F(\epsilon_{X^*}\delta_X) = \epsilon_{X^*}\delta_X,
\]
we have 
\[
d'_L = \lambda^{-1}\mu d_L, \quad
d'_R = \lambda\mu^{-1} d_R
\]
in the case $F(X) = X$ and 
\[
d'_L = \lambda \mu^{-1} d_R, \quad 
d'_R = \lambda^{-1}\mu d_L
\]
in the case $F(X) = X^*$. 

Conversely, these relations ensure the existence of 
an equivalence functor $F$ by the generating property 
of basic arcs (Lemma~2.1). 
In this way, we have proved the following: 

\begin{Proposition}
We have 
\[
\cO_{d_L,d_R} \cong \cO_{d'_L,d'_R} 
\iff d_L d_R = d'_L d'_R.
\]

In particular, if we write $\cO_d = \cO_{d,d}$, 
$\cO_{d_L,d_R} \cong \cO_d$ for the choice 
$d = \pm \sqrt{d_L d_R}$,  
with the condition 
$\cO_d \cong \cO_{d'}$ equivalent to $d = \pm d'$. 
Moreover self-equivalences of $\cO_d$ are given by 
\[
F(\epsilon_X) = \lambda \epsilon_X, 
\qquad 
F(\epsilon_{X^*}) = \lambda \epsilon_{X^*}
\]
or
\[
F(\epsilon_X) = \lambda \epsilon_{X^*}, 
\qquad 
F(\epsilon_{X^*}) = \lambda \epsilon_X
\]
with $\lambda \in \C^\times$. 
\end{Proposition}

In what follows, we shall concentrate on the category 
$\cO_d = \cO_{d,d}$
withought loss of generality. 
As in the case of Temperley-Lieb category, we then have 
the natural operation of duality on the category $\cO_d$: 
if we introduce the antimultiplicative involution $*$ on 
the monoid $\N * \N$ by switching the role of two submonoids 
$\N$, then $X^{w^*}$ gives a dual object of $X^w$ for 
$w \in \N * \N$ with respect to the obvious pairings and 
copairings by multiple arcs. The associated operation of transposed 
maps is given by rotating diagrams by an angle of $\pi$. 

\section{Semisimplicity Analysis}
In the tensor category $\cO_d$, the semisimplicity analysis, 
together with the Jones-Wenzl recursive formula, 
works for objects of alternating tensor products
$X\otimes X^*\otimes X\otimes \dots$ exactly as 
in the case of Temperley-Lieb categories (\cite{CDA}). 
In particular, all of the alternating algebras 
\[
\End(X\otimes X^*\otimes X\otimes \dots), 
\quad 
\End(X^*\otimes X\otimes X^*\otimes \dots)
\]
are semisimple 
if and only if $d = q + q^{-1}$ with $q^2$ not a proper root of 
unity.
Moreover, 
under the assumption of semisimplicity on alternating algebras 
$\End(X\otimes X^*\otimes X\otimes \dots)$, 
we can inductively define simple objects 
$\{ X_n, Y_n \}_{n \geq 1}$ 
in the idempotent-completion $\overline{\cO}_d$ of $\cO_d$
so that 
$X_n$ and $Y_n$ are the new stuffs in 
$\overbrace{X\otimes X^*\otimes \dots}^\text{$n$-factors}$ and 
$\overbrace{X^*\otimes X\otimes \dots}^\text{$n$-factors}$ 
respectively, i.e., 
$X_n = f_n (X\otimes X^*\otimes \dots)$ 
with $f_n$ the $n$-th Jones-Wenzl idempotent and 
similarly for $Y_n$. 

To a word $w \in \N * \N$, we associate a subobject 
$X_w$ of $X^w$ in $\overline{\cO}_d$ by replacing maximal 
alternating subwords in $w$ with the highest part just defined. 
(Note that $X_w^* = X_{w^*}$ and the transposed ${}^tf_n$ is 
again a Jones-Wenzl idempotent.)

\begin{Example}
For the choice 
\[
X^w = (X)(XX^*XX^*)(X^*XX^*XX^*)(X^*)(X^*),
\]
we have $X_w = X_1 \otimes X_4\otimes Y_5\otimes Y_1\otimes Y_1$. 
\end{Example}

\begin{Lemma}
For a non-empty word $w \in \N * \N$, 
$\Hom(X_w,I) = \{ 0\}$. 
\end{Lemma}

\begin{proof}
Let $w = w_1\dots w_l$ be the factorization into 
maximal alternating parts and write 
$X_{w_j} = P_j X^{w_j}$ with $P_j$ the Jones-Wenzl projection 
to the highest part of the alternating tensor product $X^{w_j}$. 
We shall show that any diagram $D$ in 
$K_{w,\emptyset} \subset \Hom(X^w,I)$ annihilates 
$P_1\otimes \dots \otimes P_l$ by an induction on the word length 
$|w| = \sum_j |w_j|$ ($|w_j| = |w_j|_+ + |w_j|_-$).

First observe that, if $D$ contains an arc connecting two vertices 
in some $w_j$, then $P_j$ annihilates the arc (regarded as 
a pairing morphism) by the highest assumption 
(old stuffs being killed by $P_j$). 
Thus, we need to deal with $D$ having no such arcs. 
Then any vertex inside $w_1$ should be joined to a vertex in 
$w_2\dots w_l$. The right-end vertex $v_1$ of $w_1$, however, 
cannot be connected to any vertex in $w_2$ 
because of parity mismatch. So $v_1$ is joined to a vertex $v_k$ 
in $w_k$ with $k \geq 3$. 
Now write $w_k = w'_k w''_k$, where $w''_k$ is the subword of $w_k$ 
starting at the vertex $v_k$ until the end of $w_k$ and 
$w'_k$ is the complement to $w''_k$. 
We then have the factorization 
$D = D''
(1_{X^{w_1}}
\otimes D'\otimes 
1_{X^{w''}})$
with $w'' = w_k''w_{k+1}\dots w_l$ and 
$D' \in K_{w_2\dots w_{k-1}w'_k}$ (see Fig.~1). 

On the other hand, 
the Clebsh-Gordan fusion rule on alternating parts ensures 
the relation $P_k = (P'_k\otimes P''_k)P_k$ 
($P_k'$ and $P_k''$ being the Jones-Wenzl idempotents for 
$w_k'$ and $w_k''$ respectively) and hence we see 
\[
D(P_1\otimes \dots \otimes P_l) = D''
(1_{X^{w_1}}\otimes D'(P_2\otimes \dots \otimes P_{k-1}\otimes P'_k) 
\otimes 1_{X^{w''}}) 
(P_1\otimes \dots \otimes P_l), 
\]
which vanishes by the induction hypothesis  
$D'(P_2\otimes \dots \otimes P_{k-1}\otimes P'_k) = 0$.
\end{proof}

\begin{figure}[h]
\input okd32.tpc
\caption{\label{okd32}}
\end{figure}

\begin{Corollary}
We have $\Hom(X_w, X_{w'}) = \{ 0\}$ for $w \not= w' \in \N * \N$ 
and 
$\Hom(X_w,X_w) = \C 1_{X_w}$. 
\end{Corollary}

\begin{proof}
Decompose $w = w_1\dots w_l$ and $w' = w'_1\dots w'_{l'}$ as 
before. Then 
$\Hom(X_w, X_{w'}) \cong \Hom(X_w\otimes X_{(w')^*},I)$ by rigidity. 
If this vector space is non-trivial, the product 
$w_l(w'_{l'})^*$ should interact at the contact point and 
we have the decomposition of the form 
\[
X_{w_l}\otimes X_{(w'_{l'})^*} \cong 
X_{\widetilde w_1} \oplus \dots \oplus X_{\widetilde w_m}
\]
according to the Clebsh-Gordan rule. 

Thus, if the unit object $I$ does not appear inside 
$X_{w_l}\otimes X_{(w'_{l'})^*}$, we know 
\[
\Hom(X_w\otimes (X_{w'})^*,I) = 
\bigoplus_{j=1}^m 
\Hom(X_{w_1\dots w_{l-1} \widetilde w_j (w'_{l'-1})^* 
\dots (w'_1)^*},I) = \{ 0\}. 
\]
When $I$ is contained, $w_l = w'_{l'}$ and the unit object 
appears exactly once in 
$X_{w_l}\otimes X_{w'_{l'}}^*$ (say, $X_{\widetilde w_1} = I$ 
and $X_{\widetilde w_j} \not\cong I$ for $j \geq 2$) and 
we get 
\begin{align*}
\Hom(X_w\otimes (X_{w'})^*,I) &\cong 
\bigoplus_{j=1}^m 
\Hom(X_{w_1\dots w_{l-1}}\otimes X_{\widetilde w_j} 
\otimes X_{w'_1\dots w'_{l'-1}}^*,I)\\
&= \Hom(X_{w_1\dots w_{l-1}}\otimes X_{w'_1\dots w'_{l'-1}}^*,I)\\
&\cong \Hom(X_{w_1\dots w_{l-1}}, X_{w'_1\dots w'_{l'-1}}).
\end{align*}
Now the induction argument is applied to see $w = w'$ and 
\[
\End(X_{w_1\dots w_l}) \cong \End(X_{w_1\dots w_{l-1}}) 
\cong \C.
\]
\end{proof}

\begin{Proposition}
The tensor category $\cO_d$ is semisimple if and only if 
$d = q + q^{-1}$ with $q^2$ not a proper root of unit. 
In this case, $\{ X_w \}_{w \in \N *\N}$ gives 
a representative set of simple objects and the fusion rule is 
given by the following recipe: 
Let $w = w_1\dots w_l$ and $w' = w'_1\dots w'_{l'}$ be 
the decompositions into maximally alternating parts.
\begin{enumerate}
\item 
If $w_lw'_1$ matches in the parity, 
the tensor product $X_{w_l}\otimes X_{w_1'}$ decomposes according to 
the Clebsh-Gordan rule, otherwise $X_w\otimes X_{w'} = X_{ww'}$ 
remains simple. 
\item
Let $X_{w_l}\otimes X_{w'_1} \cong 
X_{u_1} \oplus \dots \oplus X_{u_k}$. If $I \not\cong U_j$ for 
$1 \leq j \leq k$, 
\[
X_w\otimes X_{w'} \cong 
\bigoplus_{j=1}^k X_{\widetilde w_j}
\]
with $\widetilde w_j 
= w_1\dots w_{l-1}u_jw'_2\dots w'_{l'} \in \N * \N$ 
gives an irreducible decompositionn. 

If $X_{u_1} \not\cong I$ 
(and hence $X_{u_j} \not\cong I$ for $j \geq 2$), 
\[
X_w\otimes X_{w'} \cong 
(X_{w_1\dots w_{l-1}}\otimes X_{w'_2\dots w'_{l'}})
\oplus X_{\widetilde w_2} \oplus \dots \oplus 
X_{\widetilde w_k}
\]
and we are reduced to the decompotion of 
$X_{w_1\dots w_{l-1}} \otimes X_{w'_2\dots w'_{l'}}$.
\end{enumerate}
\end{Proposition}

\begin{Remark}
If we restrict ourselves to the unitary (i.e., C*-) case, 
then the above formula produces the fusion rule in 
\cite[Theorem~1]{Ba2} via Proposition~5.3 below. 
\end{Remark}

\section{Positivity Condition}
We shall here investigate possible *-structures on $\cO_d$. 
Assume that there is a (compatible) *-structure on $\cO_d$. 
Then we should have 
\[
(\epsilon_X)^* = c_X \delta_{X^*}, 
\quad 
(\epsilon_{X^*})^* = c_{X^*} \delta_X
\]
with $c_X, c_{X^*} \in \C^\times$.


To preserve the rigidity, we should have 
\[
1_X = 1_X^* = 
\Bigl( 
(1_X\otimes \epsilon_{X^*})(\delta_{X^*}\otimes 1_X)
\Bigr)^* = c_{X^*}(\delta_{X^*}^*\otimes 1_X) 
(1_X\otimes \delta_X)
\]
and 
\[
1_{X^*} = (1_{X^*})^* = 
\Bigl( 
(1_{X^*}\otimes \epsilon_X)(\delta_X\otimes 1_{X^*})
\Bigr)^* = c_X(\delta_X^*\otimes 1_{X^*}) 
(1_{X^*}\otimes \delta_{X^*}),
\]
which are equivalent to 
\[
(\delta_X)^* = c_X^{-1} \epsilon_{X^*}, 
\quad 
(\delta_{X^*})^* = c_{X^*}^{-1} \epsilon_X.
\]
To preserve trace (loop) values, we should have 
\begin{gather*}
\overline d = (\epsilon_X\delta_{X^*})^* 
= c_X c_{X^*}^{-1} \epsilon_X \delta_{X^*} = 
c_X c_{X^*}^{-1} d,\\
\overline d = (\epsilon_{X^*}\delta_X)^* 
= c_X^{-1} c_{X^*} \epsilon_{X^*} \delta_X = 
c_X^{-1} c_{X^*} d,
\end{gather*}
i.e., 
\[
c_X = \pm c_{X^*}\quad \text{and}\quad d = \pm \overline d.
\]
With these conditions satisfied, the universality on rigidity  
(Lemma~2.1) allows us to extend the operation to the whole hom-sets 
by antilinearity and antimultiplicativity. The obtained map is 
then involutive if and only if 
\[
\epsilon_X = (\epsilon_X)^{**} = 
\overline{c_X} \delta_{X^*}^* = \frac{\overline{c_X}}{c_{X^*}} 
\epsilon_X,\qquad
\epsilon_{X^*} = (\epsilon_{X^*})^{**} = 
\overline{c_{X^*}} \delta_X^* = \frac{\overline{c_{X^*}}}{c_X} 
\epsilon_{X^*}, 
\]
i.e., $\overline{c_X} = c_{X^*}$. 

\begin{Proposition}
The tensor category $\cO_d$ admits a compatible *-structure 
if and only if $d^2 \in \R$, i.e., $d = \pm \overline{d}$. 

If this is the case, *-structures are parametrized by 
$c \in \C^\times$ satisfying $\overline{cd} = cd$ 
with the associated *-operation given by 
\[
D^* = c^{\sharp(D)} 
\left( 
\frac{\overline d}{d}
\right)^{l(D)} D'
\]
for an oriented diagram $D$, where $D'$ is the orientation reversion 
of the reflection of $D$ (see Fig.~2), 
\[
\sharp(D) = \sharp\{ 
\text{$\epsilon_X$'s and $\epsilon_{X^*}$'s inside $D$}\}
- \sharp\{ 
\text{$\delta_X$'s and $\delta_{X^*}$'s inside $D$}\},
\]
and 
\[
l(D) = \sharp\{ 
\text{$\epsilon_{X^*}$'s inside $D$} \}
- \sharp\{ 
\text{$\delta_{X^*}$'s inside $D$} \}.
\]

Moreover two *-structures are equivalent if and only if 
$c = \lambda c'$ or $\overline c = \lambda c'$ for some 
$\lambda >0$. 
\end{Proposition}

\begin{proof}
Let $c, c' \in \C^\times$ be parameters of *-structures 
with the associated *-structures on $\cO_d$ 
denoted by $*$ and $\star$ 
respectively. Since monoidal automorphisms of $\cO_d$ are of the 
form 
\[
F(\epsilon_X) = \lambda \epsilon_X, 
\quad
F(\epsilon_{X^*}) = \lambda \epsilon_{X^*}
\quad 
\text{or}
\quad 
G(\epsilon_X) = \lambda \epsilon_{X^*}, 
\quad
G(\epsilon_{X^*}) = \lambda \epsilon_X,
\]
the condition $F(D^*) = F(D)^\star$ 
(resp.~$G(D^*) = G(D)^\star$) is equivalent to $c = |\lambda|^2 c'$ 
(resp.~$\overline c = |\lambda|^2 c'$). 
\end{proof}

\begin{figure}[h]
\hspace{1cm}
\input okd41.tpc
\caption{\label{okd41}}
\vspace{5mm}
\end{figure}

\begin{Theorem}
The tensor category $\cO_d$ admits a compatible C*-structure 
if and only if $d \in \R^\times$ and $d^2 \geq 4$. 

If this is the case, C*-structure is unique up to 
monoidal equivalences 
and is given by 
\[
D^* = \left(
\frac{d}{|d|} 
\right)^{\sharp(D)} D'
\]
for a diagram $D$. 
\end{Theorem}

\begin{proof}
We argue as in the case of Temperley-Lieb category: Since 
idempotents 
\[
e_1 = \frac{1}{d} (\delta_{X^*}\epsilon_X)\otimes 1_X, 
e_2 = \frac{1}{d} 1_X\otimes (\delta_X\epsilon_{X^*}) 
\]
are related to the central decomposition of 
$\End(X\otimes X^*)$ and $\End(X^*\otimes X)$, we have 
$e_1^* = e_1$ and $e_2^* = e_2$, which are equivalent to 
$\overline c/c = \overline d/d$ (automatically satisfied). 
From the relation 
$e_1e_2e_1 = d^{-2} e_1$, we see $d^2 >0 \iff d \in \R$ and hence 
$c \in \R^\times$ as well. Moreover $\epsilon_X\epsilon_X^* = cd>0$ 
implies the condition $c/|c| = d/|d|$. 
Thus, up to equivalences, C*-structure (if any exists) is unique 
and given by the above formula for the choice $|c|=1$. 

By a standard argument based on Jones-Wenzl formula 
(cf.~\cite{CDA}), we can derive 
the condition $d^2 \geq 4$ from the positivity of the *-structure. 
On the other hand, if this is assumed, the Temperley-Lieb category 
$\cK_d$ admits the (unique) C*-structure. 
Since $\cO_d$ is a subcategory of $\cK_d$, the *-structure of 
$\cO_d$ meets the positivity and we are done.  
\end{proof}

\section{Fiber Functors}
From the universality property, 
any fiber functor $\Phi: \cO_d \to \cV ec$ is determined by 
non-degenerate bilinear forms $\Phi(\epsilon_X)$ and 
$\Phi(\epsilon_{X^*})$. 
Set $V = \Phi(X)$, $W = \Phi(X^*)$, 
$E = \Phi(\epsilon_X): V\otimes W \to \C$ and 
$F = \Phi(\epsilon_{X^*}): W\otimes V \to \C$. 
Let $\delta_E \in W\otimes V$ and 
$\delta_F \in V\otimes W$ be the associated vector: 
Given a basis $\{ v_i\}$ of $V$ and $\{ w_j\}$ of $W$, 
\[
\delta_E = \sum c_{ji} w_j\otimes v_i, 
\quad 
\delta_F = \sum d_{ij} v_i\otimes w_j,
\]
with $\{ c_{ji}\}$ and $\{ d_{ij} \}$ defined by 
\[
\sum_j E(v_k\otimes w_j)c_{ji} = \delta_{ki}, 
\quad 
\sum_i F(w_k\otimes v_i) d_{ij} = \delta_{kj}, 
\]
give the covectors, whence they satisfy 
\[
F(\delta_E) = \sum_{ij} c_{ji} F(w_j\otimes v_i) = d 
= E(\delta_F) = \sum_{ij} d_{ij} E(v_i\otimes w_j).
\]

Conversely, given invertible matrices 
$A = (a_{ij})$ with $a_{ij} = E(v_i\otimes w_j)$ and 
$B = (b_{ji})$ with $b_{ji} = F(w_j\otimes v_i)$ satisfying 
\[
\text{trace}({}^tB A^{-1}) = d = 
\text{trace}({}^tA B^{-1}),
\]
we obtain a fiber functor. 

Two fiber functors $\Phi$ and $\Phi'$ are naturally equivalent 
if and only if we can find isomorphisms of vector spaces 
\[
\Phi(X) \to \Phi'(X) = V', 
\quad 
\Phi(X^*) \to \Phi'(X^*) = W'
\]
such that the diagrams 
\[
\begin{CD}
V \otimes W @>>> V'\otimes W'\\
@V{E}VV @VV{E'}V\\
\C @= \C
\end{CD}\ ,
\qquad 
\begin{CD}
W \otimes V @>>> W'\otimes V'\\
@V{F}VV @VV{F'}V\\
\C @= \C
\end{CD}
\]
commute. 
This means that $A$ and $B$ can be replaced by 
\[
{}^tT A S, 
\qquad 
{}^tS BT
\]
with $S$ and $T$ invertible square matrices. 
Thus we may choose $A = 1$ and there remains the gauge freedom 
of ${}^tTS = 1$, i.e., 
\[
B \mapsto T^{-1} BT.
\]

\begin{Theorem}
Fiber functors on $\cO_d$ is completely classified by similarity orbits in 
\[
\{ 
B \in \text{GL}(n,\C); \text{trace}(B) = d = 
\text{trace}(B^{-1})\}.
\]
\end{Theorem}

For the description of the associated algebraic group (or 
a Hopf algebra $H$), 
we restore the matrix $A$ and choose bases $\{ v_i\}$ in $V$ and 
$\{ w_j \}$ in $W$. Then $H$ is generated, as an algebra, 
by $2n^2$ elements 
$\{ v_{ij} \}$ and $\{ w_{kl}\}$ with 
the defining relations   
\begin{gather*}
\sum_{i,k} a_{ik} v_{ij} w_{kl} = a_{jl} 1, 
\quad 
\sum_{j,i} c_{ji} w_{kj} v_{li} = c_{kl} 1,\\
\sum_{i,k} b_{ik} w_{ij} v_{kl} = b_{jl} 1, 
\quad 
\sum_{j,i} d_{ji} v_{kj} w_{li} = d_{kl} 1
\end{gather*}
corresponding to basic morphisms 
$\begin{CD} V\otimes W @>{E}>> \C @>{\delta_E}>> W\otimes V
\end{CD}$ and 
$\begin{CD} W\otimes V @>{F}>> \C @>{\delta_F}>> V\otimes W
\end{CD}$. 
The comultiplication is given by 
\[
\Delta(v_{ij}) = \sum_k v_{ik}\otimes v_{kj}, 
\quad 
\Delta(w_{kl}) = \sum_i w_{ki}\otimes w_{il}.
\]

Next let $\Phi: \cO_d \to \cH ilb$ be a unitary fiber functor.
Here $\cO_d$ ($d^2 \geq 4$) is furnished with the *-structure 
$D^* = (d/|d|)^{\sharp(D)} D'$ and $\cH ilb$ denotes 
the C*-tensor category of finite-dimensional Hilbert spaces. 
The condition $\Phi(D)^* = \Phi(D^*)$ for any $D$, i.e.,  
\[
\Phi(\epsilon_X)^* = \frac{d}{|d|}\Phi(\delta_{X^*}), 
\quad 
\Phi(\epsilon_{X^*})^* = \frac{d}{|d|} \Phi(\delta_X)
\]
is equivalent to the relation 
\[
\overline B = \frac{d}{|d|} A^{-1}
\]
with the trace-value condition given by 
\[
\text{trace}(BB^*) = |d| = 
\text{trace}((BB^*)^{-1}). 
\]
As in the case of non-unitary fiber functors, 
the gauge freedom is given by 
\[
A \mapsto {}^tTAS, 
\qquad 
B \mapsto {}^tS BT
\]
with $S$, $T$ unitary matrices. 
Note that the condition $\overline A = (d/|d|) B^{-1}$ as well as 
the trace-value relation is preserved under this transformation. 
Thus the orbit space for invertible matrices $B$ is identified with 
\[
U(n)\backslash \text{GL}(n,\C)/U(n) 
\cong \{ \text{positive definite matrices} \}/ 
\text{unitary similarity}.
\]

\begin{Theorem}(cf.~\cite[Theorem~1.1]{Wa}, \cite[Theorem~6.2]{BDV})
Unitary fiber functors on the C*-tensor category $\cO_d$ 
are parametrized by positive eigenvalue lists 
$\{ \mu_1, \dots, \mu_n\}$ satisfying 
\[
\sum_j \mu_j^2 = |d| = \sum_j \mu_j^{-2}.
\]
\end{Theorem}

We shall identify the associated compact quantum group 
with the universal quantum group of unitary type 
in \cite{Ba, D-W}. 

For this, we use the notation and the construction 
in \cite{FTL}. Choose orthonormal bases $\xi = \{ v_i\}$ for 
$V$ and $\eta = \{ w_k\}$ for $W$ with 
the basis $\xi^* = \{ v^*_j\}$ of $W$ defined 
by the relation 
$E(v_j^*\otimes v_i) = \delta_{ji}$, i.e., 
$\sum_k (w_k|v^*_i) b_{kj} = \delta_{ij}$. 
Let $\{ u^\xi_{ij} \}$ and $\{ u^\eta_{kl} \}$ be the associated 
generators.  
Then, from the covariance condition,
\[
\sum_l (w_l|v^*_i) u^\eta_{kl} = 
\sum_j (w_k|v^*_j) u^{\xi^*}_{ji}
\]
and the *-relation $(u^\xi_{ij})^* = u^{\xi^*}_{ij}$, we see that 
\[
u^\eta ({}^tB^{-1}) = ({}^tB^{-1}) u^{\xi^*}
\]
in the matrix notation for $u^\eta$ and $u^{\xi^*}$. 
Thus the associated Hopf *-algebra $H$ is generated by 
$u = (u^\xi_{ij})$, which constitutes a unitary matrix 
by the orthonormality of the basis $\xi = \{ v_j\}$. 

From the universality on rigidity, the defining relations 
are given by the covariance conditions for the morphisms
$\Phi(\epsilon_X)$, $\Phi(\epsilon_{X^*})$, 
$\Phi(\delta_X)$ and $\Phi(\delta_{X^*})$:
\begin{align*}
\Phi(\epsilon_X)(v_i\otimes w_k) 1_H &= 
\sum \Phi(\epsilon_X)(v_j\otimes w_l) 
u^{\xi\otimes \eta}_{j,l;i,k},\\ 
\Phi(\epsilon_{X^*})(w_k\otimes v_i) 1_H &= 
\sum \Phi(\epsilon_{X^*})(w_l\otimes v_j) 
u^{\eta\otimes\xi}_{l,j;k,i},\\
\sum u^{\eta\otimes \xi}_{k,i;l,j} 
(w_l\otimes v_j| \Phi(\delta_X)) &= 
(w_k\otimes v_i| \Phi(\delta_X)) 1_H,\\
\sum u_{i,k;j,l}^{\xi\otimes \eta} 
(v_j\otimes w_l|\Phi(\delta_{X^*})) &= 
(v_i\otimes w_k| \Phi(\delta_{X^*})) 1_H. 
\end{align*}
From the multiplication relations 
$u^{\xi\otimes \eta}_{i,k;j,l} = u^\xi_{ij} u^\eta_{kl}$, 
$u^{\eta\otimes \xi}_{k,i;l,j} = u^\eta_{kl} u^\xi_{ij}$, 
the condition is equivalent to 
\begin{align*}
A &= {}^t u^\xi A u^\eta,\\
B &= {}^t u^\eta B u^\xi,\\
u^\eta A^{-1} {}^t u^\xi &= A^{-1},\\
u^\xi B^{-1} {}^t u^\eta &= B^{-1}
\end{align*}
and therefore to the conditions $uu^* = 1 = u^*u$, 
$\overline u ({}^tBA^{-1}) {}^t u = {}^tB A^{-1}$ and  
${}^t u A {}^tB^{-1} \overline u = A {}^tB^{-1}$.
Since $A^{-1} = (d/|d|) \overline B$, the conditions are further 
equivalent to 
\[
uu^* = 1 = u^* u,
\quad
\overline u (AA^*)^{-1} {}^t u = (AA^*)^{-1}, 
\quad 
{}^t u AA^* \overline u = AA^*.
\]
Thus the associated compact quantum group is identified with 
$A_u(AA^*)$ in \cite{D-W}. 

\begin{Proposition}
The compact quantum group associated to the unitary fiber functor 
is naturally isomorphic to $A_u(AA^*)$ in \cite{D-W}, 
i.e., $A_u(A^*)$ in \cite{Ba}. 
\end{Proposition}

\begin{Remark}
In \cite{BDV}, unitary fiber functors are classified for 
representation categories of the compact quantum group $A_u(F)$ 
based on representation theory of \cite{Ba2}. 
Since these representation categories are exactly $\cO_d$ as tensor 
categories by the above proposition, we get an access to some results 
in \cite{BDV} in an elementary way.  
\end{Remark}

\appendix
\section{Faithfulness}
We shall here check the automatic faithfulness of relevant 
functors. 

\begin{Proposition}
Fiber functors on the Temperley-Lieb categories 
are faithful whenever the fundamental vector space $V$ has 
dimension two or more.

Fiber functors on $\cO_d$ are faithful if 
the fundamental vector spaces $V$ and $W$ have dimension two or more.
\end{Proposition}

\begin{proof}
Since Frobenius transforms are isomorphisms, we need to show that 
the family $\{ \Phi(D); D \in K_{2n,0} \}$ is linearly indenpendent 
for any $n \geq 1$, where 
$\Phi: \cT\cL_d \to \cV ec$ denotes a fiber functor such that 
the fundamental vector space $V = \Phi(X)$ has dimension 
two or more. 

Given $1 \leq k \leq n$, let $\{ x_1, \dots x_k \}$ be 
the first $k$-vertices (counting from the left end) 
for diagrams in $K_{2n,0}$ 
and set 
\[
\cD_{k,n} = \{ 
D \in K_{2n,0}; 
\text{there are no arcs connecting $x_i$ and $x_j$ in $D$}
\}. 
\]
Note that $\cD_{n,n} \subset \cD_{n-1,n} \subset \dots 
\subset \cD_{2,n} \subset \cD_{1,n} = K_{2n,0}$ and 
$\cD_{n,n}$ consists of a single diagram. 

Since $\dim V \geq 2$, 
we can find a vector $0 \not= v \in V$ such that 
$F(\epsilon)(v\otimes v) = 0$ 
as a solution of a quadratic equation, for which 
we shall show that the family 
\[
\cE_{k,n} = \{ \Phi(D)(v^{\otimes k}\otimes \cdot)
V^{\otimes (2n-k)} \to \C; D \in \cD_{k,n} \}
\]
is linearly independent by an induction on $(k,n)$. 
Observe that, if the family $\cD_{1,n}$ is linearly independent, 
so is the family $\{ \Phi(D); D \in K_{2n,0} \}$. 

Given $1 \leq k \leq n$, 
assume the linear independence of $\cE_{k',n'}$ for 
$1 \leq k' \leq n' < n$ and for $k' > k, n' = n$. 
We shall prove that the family $\cE_{k,n}$ is linearly independent. 
Suppose that 
\[
\sum_{D \in \cD_{k,n}} c_D \Phi(D)(v^{\otimes k}\otimes \cdot) = 0. 
\]
Performing the evalution by 
the vector $v$ one step further, we have 
\[
\sum_{D \in \cD_{k,n}} c_D \Phi(D)(v^{\otimes (k+1)}\otimes \cdot) 
= 0 
\]
in $(V^{\otimes (2n-k-1)})^*$. Since $\Phi(D)$ is killed by 
this evaluation for $D \in \cD_{k,n} \setminus \cD_{k+1,n}$, 
the summation can be restricted to $\cD_{k+1,n}$ 
and the induction hypothesis ensures $c_D = 0$ for 
$D \in \cD_{k+1,n}$. 

By the non-degeneracy of the bilinear form $\Phi(\epsilon)$, 
we can find a vector $v' \in V$ satifying 
$\Phi(\epsilon)(v\otimes v') = 1$. Evaluating by $v'$ in 
the starting equation, we then have 
\[
\sum_{D \in \cD_{k,n} \setminus \cD_{k+1,n}}
c_D \Phi(D)(v^{\otimes k}\otimes v'\otimes \cdot) = 0
\]

Since $\cD_{k,n} \setminus \cD_{k+1,n}$ consists of diagrams 
which contain the arc connecting $x_k$ and $x_{k+1}$ 
(Fig.~3), 
we have the natural bijection 
$\cD_{k,n} \setminus \cD_{k+1,n} \ni D \mapsto 
D' \in \cD_{k-1,n-1}$ and then the last equation takes the form 
\[
\sum_{D' \in \cD_{k-1,n-1}} c_D 
\Phi(D')(v^{\otimes (k-1)}\otimes \cdot) = 0.
\]
Again by the induction hypothesis, 
$c_D = 0 $ for $D \in \cD_{k,n} \setminus \cD_{k+1,n}$ and 
we are done for the Temperley-Lieb case. 

For a fiber functor $\Phi: \cO_d \to \cV ec$, we normalize 
the bilinear forms $\Phi(\epsilon_X): V\otimes W \to \C$ and 
$\Phi(\epsilon_Y): W\otimes V \to \C$ so that the associated 
matrices $A$ and $B$ are upper-triangular. 
Then, for the choice $v = (1,0, \dots, 0)$ and 
$w = (0,1,0,\dots,0)$, we see 
$\Phi(\epsilon_X)(v\otimes w) = \Phi(\epsilon_Y)(w\otimes v) = 0$. 

Now, given an object $X^\omega$ of $\cO_d$, we can take the tensor 
product of $v$ and $w$ according to the arrangement of $X$ and 
$Y = X^*$ in $X^\omega$ up to the $k$-th factor. 
By using these as probing vectors, we can repeat the above argument
to conclude the independence of the family 
$\{ \Phi(D); D \in K_{\omega, \emptyset}\}$. 

\begin{figure}[h]
\hspace{1cm}
\input okdaa.tpc
\caption{\label{okdaa}}
\vspace{5mm}
\end{figure}

\end{proof}


\begin{thebibliography}{8}
\bibitem{Ba}
T.~Banica, 
Th\'eorie des repr\'esentations du groupe quantique compact libre 
$O(n)$, 
\textit{C.~R.~Acad.~Sci.~Paris}, 
322(1996), 241--244. 
\bibitem{Ba2}
\underline{\phantom{T.~Banica}}, 
Le groupe quantique compact libre $U(n)$, 
\textit{Commun.~Math.Phys.}, 
190(1997), 143--172. 
\bibitem{BDV}
J.~Bichon, A.~De Rijdt and S.~Vaes, 
Ergodic coactions with large multiplicity and monoidal
 equivalence of quantum groups, 
\textit{Commun.~Math.Phys.}, 
262(2006), 703--728. 
\bibitem{D-W}
A.~van Daele and S.Z.~Wang, 
Universal quantum groups, 
\textit{Intern.~J.~Math.}, 7(1996), 255--264. 
\bibitem{Kau}
L.H.~Kauffman, State models and the Jones polynomial, 
{\it Topology}, 
26(1987), 395--407.
\bibitem{Kau2}
\underline{\phantom{Kauffman}}, 
An invariant of regular isotopy, 
\textit{Trans.~Amer.~Math.~Soc.}, 
318(1990), 417--471.
\bibitem{Wa}
S.~Wang, Structure and isomorphism classification of 
compact quantum groups $A_u(Q)$ and $B_u(Q)$, 
\textit{J.~Operator Theory}, 48(2002), 573--583.
\bibitem{CDA}
S.~Yamagami, 
A categorical and diagrammatical approach to Temperley-Lieb 
algebras, 
math.QA/0405267.
\bibitem{FTL}
\underline{\phantom{S.~Yamagami}}, 
Fiber functors on Temperley-Lieb categories, 
math.QA/0405517.
\end{thebibliography}
\end{document}